\documentclass[12pt]{amsart}
\usepackage{graphicx}
\begin{document}
\newtheorem{theorem}{Theorem}
\newtheorem{prop}[theorem]{Proposition}
\newtheorem{lemma}[theorem]{Lemma}
\newtheorem{claim}[theorem]{Claim}
\newtheorem{cor}[theorem]{Corollary}
\newtheorem{defin}[theorem]{Definition}
\newtheorem{defins}[theorem]{Definitions}
\newtheorem{example}[theorem]{Example}
\newtheorem{xca}[theorem]{Exercise}
\newcommand{\map}{\mbox{$\rightarrow$}}
\newcommand{\aaa}{\mbox{$\alpha$}}
\newcommand{\Aaa}{\mbox{$\mathcal A$}}
\newcommand{\bbb}{\mbox{$\beta$}}
\newcommand{\ccc}{\mbox{$\mathcal C$}}
\newcommand{\ddd}{\mbox{$\delta$}}
\newcommand{\Ddd}{\mbox{$\Delta$}}
\newcommand{\eee}{\mbox{$\epsilon$}}
\newcommand{\Fff}{\mbox{$\mathcal F$}}
\newcommand{\Ggg}{\mbox{$\Gamma$}}
\newcommand{\ggg}{\mbox{$\gamma$}}
\newcommand{\Hhh}{\mbox{$\mathcal H$}}
\newcommand{\kkk}{\mbox{$\kappa$}}
\newcommand{\lll}{\mbox{$\lambda$}}
\newcommand{\Lll}{\mbox{$\Lambda$}}
\newcommand{\ooo}{\mbox{$\omega$}}
\newcommand{\Rr}{\mbox{$\mathbb{R}$}}
\newcommand{\Rrr}{\mbox{$\mathbb{R}^{2}$}}
\newcommand{\rrr}{\mbox{$\rho$}}
\newcommand{\sss}{\mbox{$\sigma$}}
\newcommand{\Sss}{\mbox{$\mathcal S$}}
\newcommand{\Ss}{\mbox{$\Sigma$}}
\newcommand{\Th}{\mbox{$\Theta$}}
\newcommand{\ttt}{\mbox{$\tau$}}
\newcommand{\bdd}{\mbox{$\partial$}}
\newcommand{\zzz}{\mbox{$\zeta$}}
\newcommand{\Zz}{\mbox{$\mathbb{Z}$}}
\newcommand{\Zzz}{\mbox{$\mathbb{Z}^{2}$}}

\newcommand{\inter}{\mbox{${\rm int}$}}

\title[] {Automorphisms of the 3-sphere that preserve a genus two 
Heegaard splitting}

\author{Martin Scharlemann}
\address{\hskip-\parindent
         Mathematics Department\\
         University of California\\
         Santa Barbara, CA 93106\\
         USA}
\email{mgscharl@math.ucsb.edu}

\date{\today}
\thanks{Research supported in part by an NSF grant.}

\begin{abstract}
An updated proof of a 1933 theorem of Goeritz, exhibiting a
finite set of generators for the group of automorphisms of $S^{3}$
that preserve a genus two Heegaard splitting.  The group is analyzed 
via its action on a certain connected $2$-complex.
\end{abstract}

\maketitle

\section{Introduction}

In 1933 Goeritz \cite{Go} described a set of automorphisms of the
standard unknotted genus two handlebody in $S^{3}$, each of which
extends to all of $S^{3}$.  He further observed that any such
automorphism is a product of elements of this finite set.  Stated
somewhat differently, Goeritz identified a finite set of generators
for the group $\Hhh$, defined as isotopy classes of
orientation-preserving homeomorphisms of the $3$-sphere that leave a
genus two Heegaard splitting invariant.  Goeritz' theorem was
generalized to Heegaard splittings of aribtrarily high genus by Powell
\cite{Po}, but the proof contains a serious gap.\footnote{On p. 
210, Case 2 the argument requires that, among the chambers into which
$\phi^{-1}(s_k)$ divides the handle, there are two adjacent ones that
each contain pieces of $G_{k}^{h}$.  There is no apparent reason why
this should be true.}  So a foundational question remains unresolved:
Is the group of automorphisms of the standard genus $g$ Heegaard
splitting of $S^{3}$ finitely generated and, if so, what is a natural
set of generators.  The finite set of elements that Powell proposes
as generators remains a very plausible set.

Since the gap in Powell's proof has escaped attention for 25 years,
Goeritz' original theorem might itself be worth a second look.  In
addition, his argument is difficult for the modern reader to follow,
is published in a fairly inaccessible journal and is a bit
old-fashioned in its outlook.  In view of the use that has been made
of it in recent work on tunnel number one knots (cf \cite{ST},
\cite{Sc}) it seems worthwhile to present an updated proof, in hopes
also that it might be relevant to the open analogous problem for
Heegaard splittings of higher genus.

The purpose of this note is to present such a proof, one influenced by
the idea of thin position.  One way to describe the outcome of this
investigation is this: there is a natural $2$-complex $\Ggg$ (which
deformation retracts to a graph) on which $\Hhh$ acts transitively. 
One can write down an explicit finite presentation for the stabilizer
$\Hhh_{P}$ of a vertex $v_{P} \in \Ggg$ and observe that the
stabilizer acts transitively on the edges of $\Ggg$ incident to
$v_{P}$.  In particular, if we add to $\Hhh_{P}$ any element $\ddd$ of
$\Hhh$ that takes $v_{P}$ to some adjacent vertex then the subgroup
generated by $\Hhh_{P}$ and $\ddd$ is exactly the subgroup that
preserves the component in which $v_{P}$ lies.  This in fact is all of
$\Hhh$, because it turns out that $\Ggg$ is connected.  The proof that
$\Ggg$ is connected can be viewed as the core argument in this paper.

\section{The complex $\Ggg$ and its vertex stabilizers}

We outline the general setting, referring the reader to \cite[Section
1]{Po} for details.  Let $V$ denote the standard unknotted genus-two
handlebody in $S^{3}$, with closed complement $W$ also a genus two
handlebody.  Let $\Hhh$ denote the group of orientation-preserving
homeomorphisms of $S^{3}$ that preserve $V$.  Regard two as equivalent
if there is an isotopy from one to the other via isotopies that
preserve $V$.  Any orientation preserving homeomorphism of $S^{3}$ is
isotopic to the identity, so an element $h: (S^{3}, V) \map (S^{3},
V)$ of $\Hhh$ is isotopic, as a homeomorphism of $S^{3}$, to the
identity.  This gives an alternate view of $\Hhh$: an element of
$\Hhh$ corresponds to an isotopy of $S^{3}$ from the identity to a
homeomorphism that preserves $V$ setwise.

For $T = \bdd V = \bdd W$, $S^{3} = V \cup_{T} W$ is a genus two
Heegaard splitting of $S^{3}$.  In the language of Heegaard
splittings, a {\em reducing sphere} $P \subset S^{3}$ is a sphere that
intersects $T$ transversally in a single essential circle and so
intersects each handlebody in a single essential disk.  Since $P$ is
separating in $S^{3}$, $P \cap T$ is a separating curve in $T$, which
we will denote $c$.  A straightforward innermost disk argument shows
that $P$ is determined up to isotopy rel $T$ by the circle $c$.

\begin{figure} [tbh]
\centering
\includegraphics[width=0.7\textwidth]{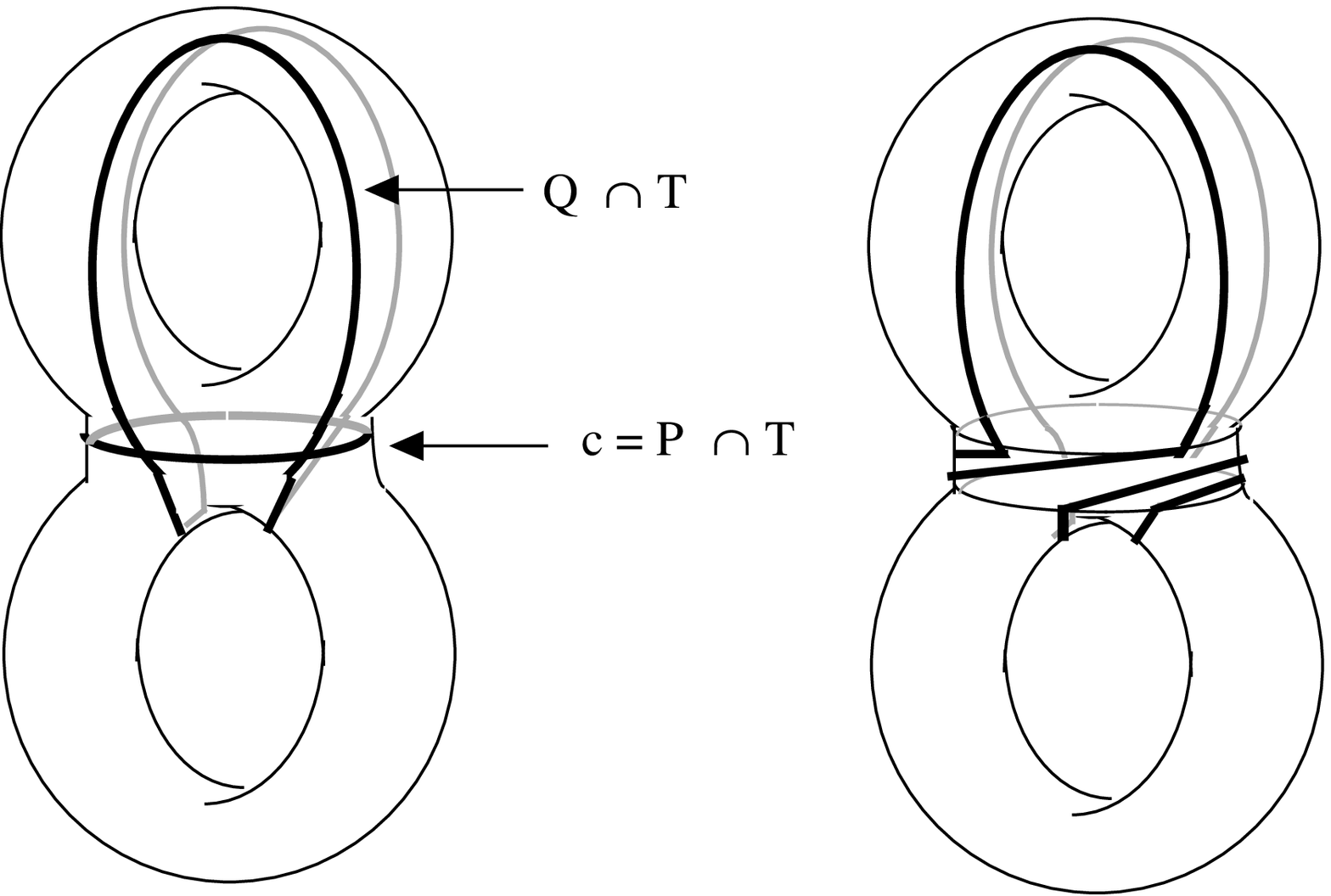}
\caption{} \label{fig:PandQ}
\end{figure}

Suppose $Q$ is another reducing sphere, with the circles $c$ and $Q
\cap T$ isotoped to intersect transversally and minimally in $T$. 
Then the number of points of intersection $|P \cap T \cap Q|$ is
denoted $P \cdot Q$.  Clearly $P \cdot Q = 0$ if and only if $P$ and
$Q$ are isotopic since the only separating essential curve in either
punctured torus component of $T - c$ is boundary parallel.  Since
reducing spheres are separating, $P \cdot Q$ is always even.  An
elementary argument (see \cite[Lemma 2.5]{ST}) shows that $P \cdot Q
\neq 2$ and in some sense characterizes (up to multiple 
half-Dehn twists about $c$) all spheres $Q$ so that $P \cdot Q = 4$. 
See Figure \ref{fig:PandQ}

\begin{figure} [tbh]
\centering
\includegraphics[width=0.4\textwidth]{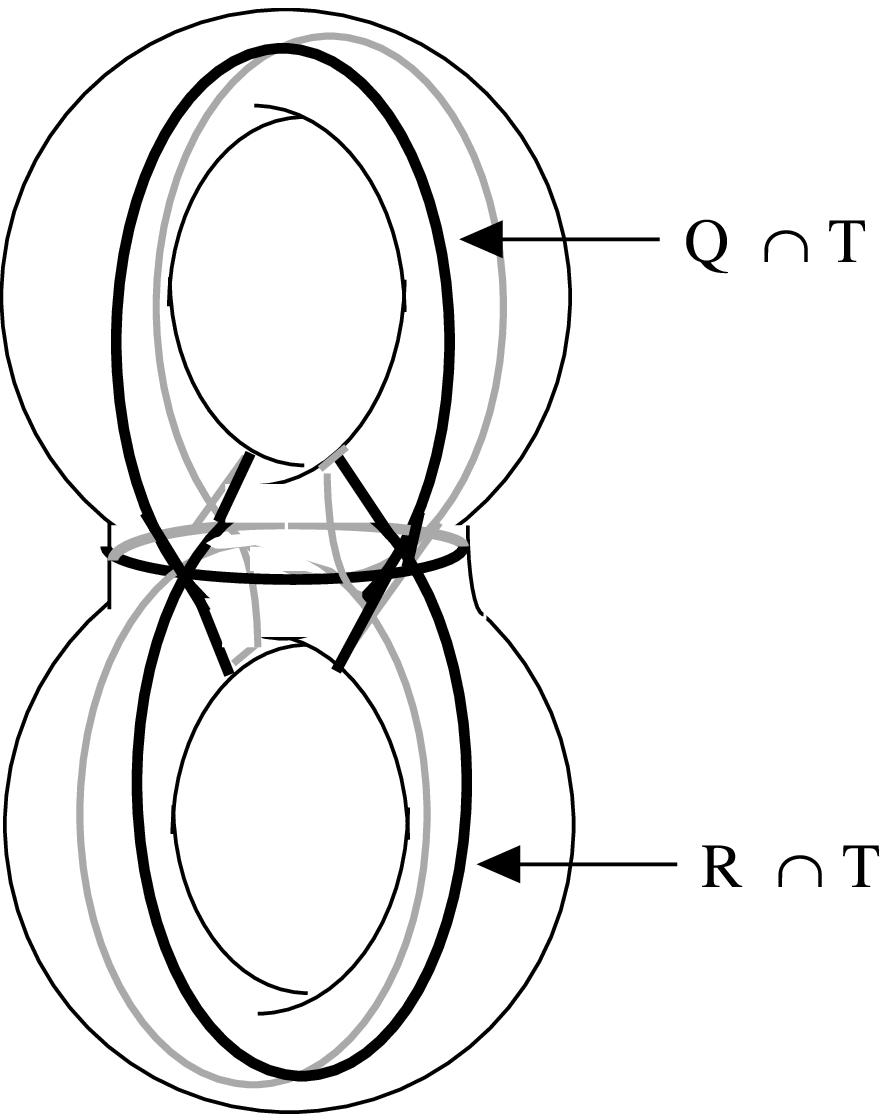}
\caption{} \label{fig:PandQandR}
\end{figure}

\begin{figure} [tbh]
\centering
\includegraphics[width=0.5\textwidth]{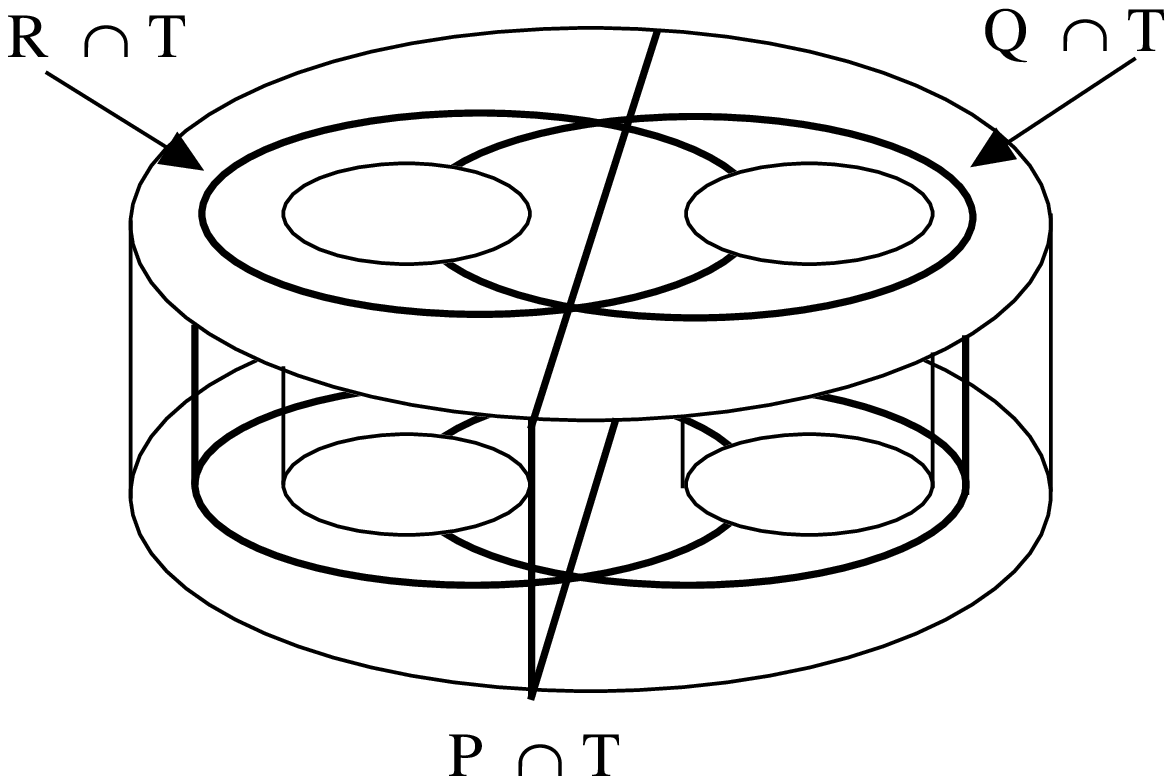}
\caption{} \label{fig:PandQandR2}
\end{figure}

This suggests a useful simplicial complex: Let $\Ggg$ be the complex
in which each vertex represents an isotopy class of reducing spheres
and a collection $P_{0}, \ldots, P_{n}$ of reducing spheres bounds an
$n$-simplex if and only if $P_{i} \cdot P_{j} = 4$ for all $ 0 \leq i
\neq j \leq n$.  In fact it follows easily from the characterization
in \cite[Lemma 2.5]{ST} that $n \leq 2$.  Figure \ref{fig:PandQandR}
illustrates a collection of three reducing spheres whose corresponding
vertices in $\Ggg$ span a $2$-complex in $\Ggg$.  (An alternate view,
in which $V$ appears as \textit{(pair of pants)} $\times I$, is shown
in Figure \ref{fig:PandQandR2}.)  Thus $\Ggg$ is a $2$-complex.

Each edge of $\Ggg$ lies on a single 2-simplex.  This is perhaps best
seen in Figure \ref{fig:PandQandR2}: The curve $P \cap T$ is uniquely
defined by the curves $Q \cap T$ and $R \cap T$ shown.  (For example,
if the curve $P \cap T$ is altered by Dehn twists around the outside
boundary of the pair of pants, it becomes a curve that is non-trivial
in $\pi_{1}(V)$, so it can't bound a disk in $V$.)  So the $2$-complex
$\Ggg$ deformation retracts naturally to a graph, in which each
$2$-simplex in $\Ggg$ is replaced by the cone on its three vertices.

A reducing sphere $P$ divides $S^{3}$ into two $3$-balls $B_{\pm}$ and
$T$ intersects each $3$-ball in a standard unknotted
punctured torus, unique up to isotopy rel boundary.  It follows that
for any two reducing spheres $P$ and $Q$ there is an orientation
preserving homeomorphism of $S^{3}$, preserving $V$ as a set, that
carries $P$ to $Q$.  Thus $\Hhh$ acts transitively on the vertices of
$\Ggg$.

\bigskip

We now explicitly give a presentation of the group that stabilizes a
vertex of $\Ggg$.  As above, let $P$ be a reducing sphere for the
Heegaard splitting $S^{3} = V \cup_{T} W$ and suppose $h: (S^{3}, V)
\map (S^{3}, V)$ is an orientation preserving homeomorphism that
leaves $P$ invariant.  That is, suppose $h$ represents an element in
$\Hhh$ that stabilizes the vertex in $\Ggg$ corresponding to $P$.

\begin{figure} [tbh]
\centering
\includegraphics[width=0.5\textwidth]{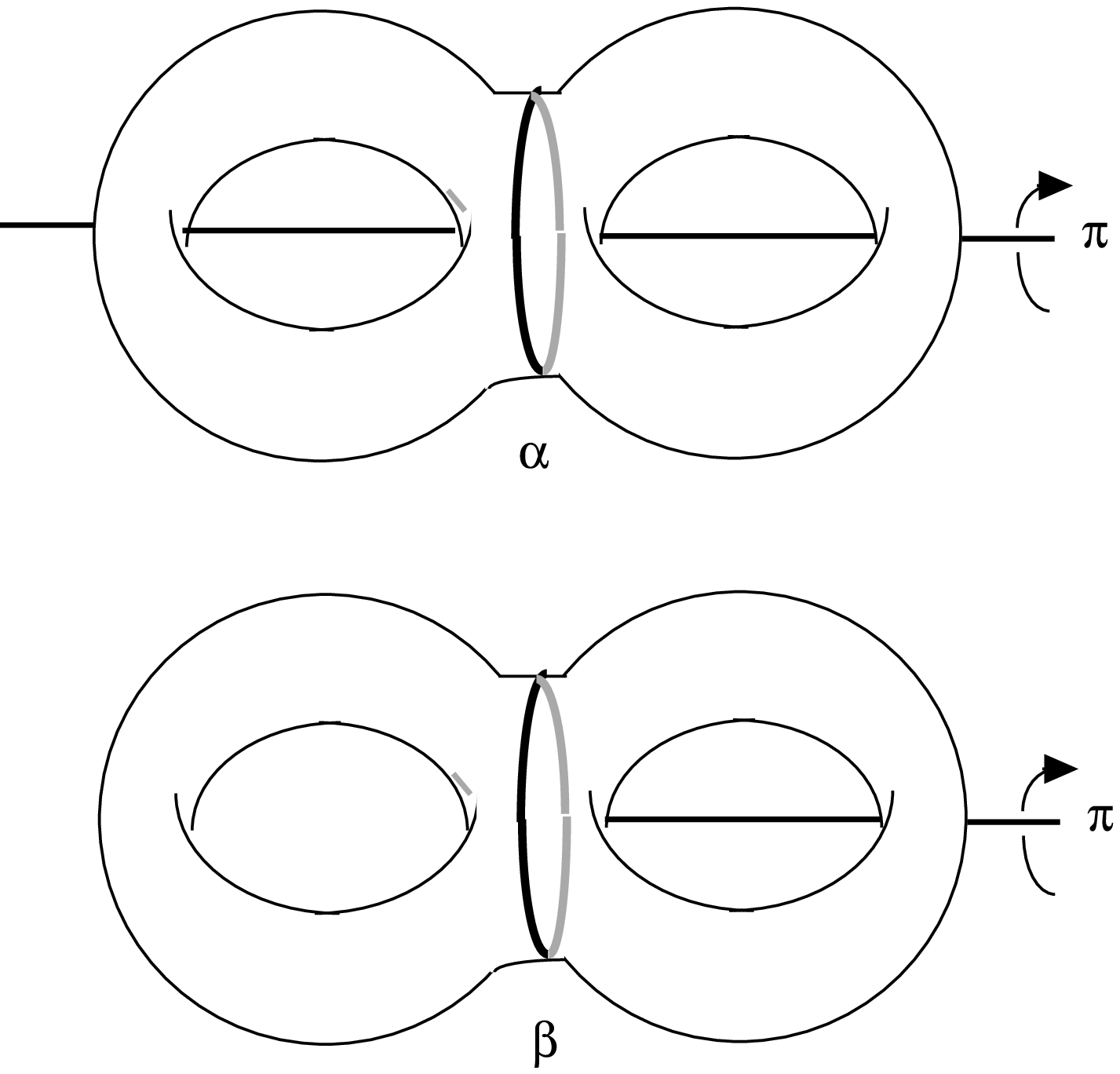}
\caption{} \label{fig:alphabeta}
\end{figure}

First assume that $h$ preserves the orientation of $P$.  Let $T_{\pm}
= T \cap B_{\pm}$ denote the two punctured torus components of $T -
P$; since $h$ preserves orientation of both $S^{3}$ and $P$ we have
$h(T_{+}) = T_{+}$ and $h(T_{-}) = T_{-}$.  Up to isotopy there is a
unique non-separating curve $\mu_{\pm} \subset T_{\pm}$ that bounds a
disk in $V$ and a unique non-separating curve $\lll_{\pm}$ that bounds
a disk in $W$ and we may choose these curves so that $\mu_{\pm} \cap
\lll_{\pm}$ is a single point.  Hence, up to equivalence in $\Hhh$, we
may with little difficulty assume that each wedge of circles
$\ggg_{\pm} = \mu_{\pm} \cup \lll_{\pm}$ is mapped to itself by $h$
and, on each $\ggg_{\pm}$, the homeomorphisms $h|\mu_{\pm}: \mu_{\pm}
\map \mu_{\pm}$ and $h|\lll_{\pm}: \lll_{\pm} \map \lll_{\pm}$ are
either simultaneously orientation preserving (in which case we can
take them both to be the identity) or simultaneously orientation
reversing (in which case we can take them each to be reflections that
preserve their intersection point).  Having identified $h$ on
$\ggg_{\pm}$ we observe that $T - (\ggg_{+} \cup \ggg_{-})$ is an
annulus $A$, and any end-preserving homeomorphism $A \map A$ is
determined up to isotopy and Dehn twists around its core by $h|\bdd
A$.  The upshot of this discussion is the following description:

\begin{lemma}: Let $\Hhh_{P}^{+}$ be the subgroup of $\Hhh$
represented by homeomorphisms that restrict to orientation-preserving
homeomorphisms of $P$.  Then $$\Hhh_{P}^{+} \cong Z_{2} + Z$$ with
generators given by the automorphisms $\aaa$ and $\bbb$ shown in
Figure \ref{fig:alphabeta}.
\end{lemma}

\begin{figure} [tbh]
\centering
\includegraphics[width=0.4\textwidth]{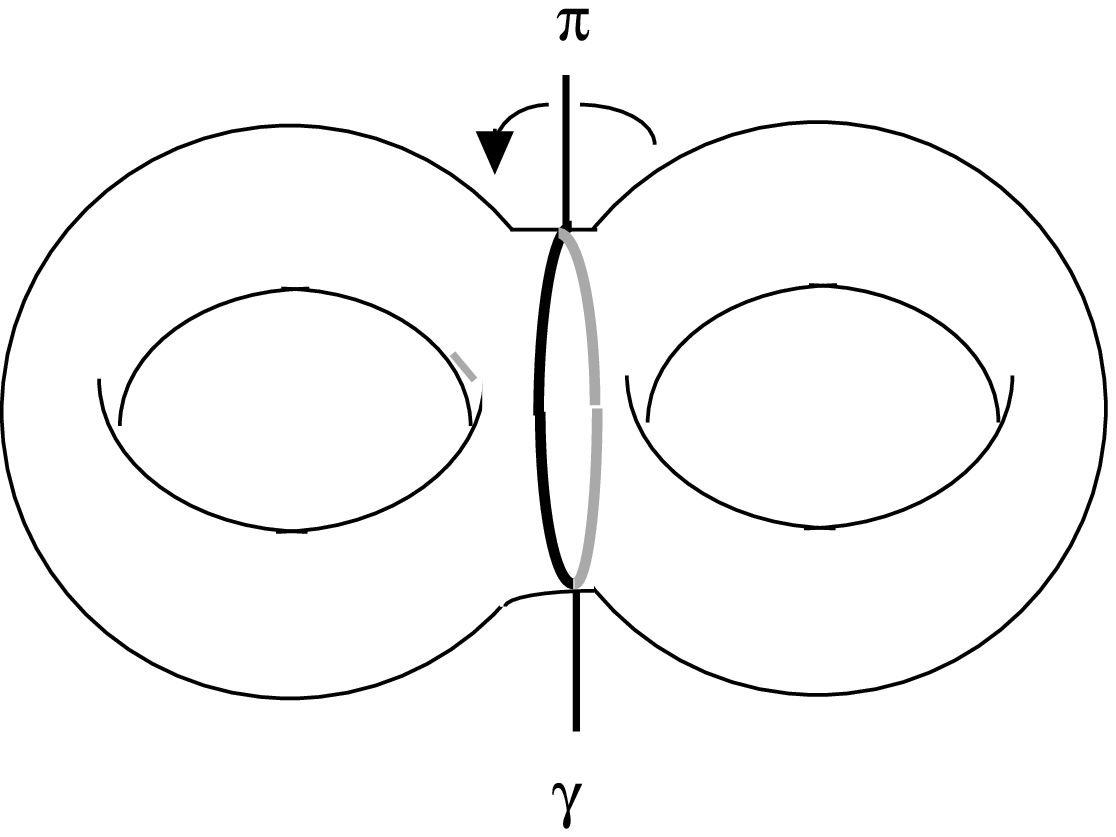}
\caption{} \label{fig:gamma}
\end{figure}

The situation is only slightly more complicated if we drop the
requirement that $h|P$ be orientation preserving since the order two
element $\ggg \in \Hhh$ shown in Figure \ref{fig:gamma} preserves $P$ but
reverses its orientation.

\begin{lemma} \label{lemma:stabilize} Let $\Hhh_{P}$ be the subgroup
of $\Hhh$ represented by homeomorphisms that preserve $P$.  Then
$\Hhh_{P}$ is an extension of $\Hhh_{P}^{+}$ by $Z_{2}$, via the
relations $\ggg \aaa \ggg = \aaa$ and $\ggg \bbb \ggg = \aaa \bbb$.
\end{lemma}

Finally, observe that if $Q$ and $Q'$ are reducing sphere so that $P
\cdot Q = 4$ and $ P \cdot Q' = 4$ then for some $n \in \Zz$, either
$\bbb^{n}$ or $\bbb^{n} \ggg$ carries $Q$ to $Q'$.  (See discussion of
Figure \ref{fig:PandQ} above.)  Interpreting this in terms of the
action of $\Hhh$ on the complex $\Ggg$ we have:

\begin{cor} \label{cor:linktransitive} Let $\Hhh_{P}$ be the subgroup
of $\Hhh$ that stabilizes the vertex $v_{P} \in \Ggg$ corresponding to
$P$.  Then $\Hhh_{P}$ is transitive on the edges of $\Ggg$ incident to
$v_{P}$.
\end{cor}

\section{Intersection of reducing spheres}

Suppose $T_{0}$ is an oriented punctured torus containing oriented
simple closed curves $\mu, \lll$ that intersect in a single point. 
For $\aaa$ an essential embedded arc in $T_{0}$ define the {\em slope}
$\sss(\aaa) \in \mathbb{Q} \cup \{\infty \}$ of the arc $\aaa$ as
follows: Orient $\aaa$ and let $p = \aaa \cdot \mu$ and $q = \aaa
\cdot \lll$ be the algebraic intersection numbers of the corresponding
homology classes.  Then $\sss(\aaa) = p/q$.  Reversing the orientation
of $\aaa$ has no effect on the slope, since it changes the sign of
both $p$ and $q$.  An alternate description of the (unsigned) slope is
this: minimize by an isotopy in $T_{0}$ the numbers $p = \aaa \cap
\mu$ and $q = \aaa \cap \lll$; then $|\sss(\aaa)| = p/q$.  If $\bbb
\subset T_{0}$ is another essential arc, with slope $r/s$ define their
distance $\Ddd(\aaa, \bbb) = |ps - qr| \in \mathbb{N}$.  It is easy to
see that if the arcs $\aaa$ and $\bbb$ are disjoint then $\Ddd(\aaa,
\bbb) \leq 1$.  Any embedded collection of arcs in $T_{+}$ constitutes
at most three parallel families of arcs, with slopes of any pair of
disjoint non-parallel arcs at a distance of one.

We now apply this terminology in the setting given above: $P$ is a
reducing sphere for $V \cup_{T} W$, the closed $3$-ball components of
$S^{3} - P$ are $B_{\pm}$, the punctured tori $T \cap B_{\pm}$ are
denoted $T_{\pm}$ and $Q$ is a reducing sphere for $V \cup_{T} W$ that
is not isotopic to $P$ and has been isotoped so as to minimize $|P
\cap Q \cap T| = P \cdot Q$.  It will be convenient to imagine $P$ as
a level sphere of a standard height function on $S^{3}$, with $B_{+}$
above $P$ and $B_{-}$ below $P$.  When we use the terms {\em above}
and {\em below} in what follows, we will be refering to such a height
function.

In each of $T_{\pm}$ there are closed non-separating curves
$\mu_{\pm}, \lll_{\pm}$ bounding respectively disks in $V$ and disks
in $W$ and for each pair, $\mu_{\pm} \cap \lll_{\pm}$ is a single
point.  We will consider the collection of arcs $Q \cap T_{\pm}$ and
their slopes with respect to $\mu_{\pm}, \lll_{\pm}$.  Fix at the
outset some orientations, e.g. orient $T$ (hence $T_{\pm}$) as $\bdd
V$ and orient $\mu_{\pm}, \lll_{\pm}$ so that the algebraic
intersection number $\mu_{\pm} \cdot \lll_{\pm} = 1$.  (The exact
choice of orientations is not critical.)

\begin{lemma} \label{lemma:outermost} There is some arc in either $Q
\cap T_{+}$ or in $Q \cap T_{-}$ that is of slope $\infty$ and another
such arc is of slope $0$.
\end{lemma}

\begin{proof} An outermost disk cut off by the disk $P \cap V$ from the
disk $Q \cap V$ is a meridian disk $D$ of the
solid torus $V \cap B_{+}$ or $V \cap B_{-}$.  Then the arc $D \cap T$
must be of slope $0$.  A symmetric argument on the disks $P \cap W, Q
\cap W$ gives an arc of slope $\infty$.
\end{proof}

\begin{lemma} \label{lemma:cross} Suppose that an arc $\aaa_{+}$ of $Q
\cap T_{+}$ has slope $\infty$ (resp.  $0$) and that there is an arc
$\aaa_{-}$ of slope $0$ (resp.  $\infty$) in $T_{-}$ that is disjoint
from $Q$.  Then there is a reducing sphere $R$ so that $P \cdot R = 4$
and $R \cdot Q < P \cdot Q$.  

The same hypothesis, but with $T_{+}$ and $T_{-}$
reversed, leads to the same conclusion.
\end{lemma}

\begin{proof} Since $\aaa_{-}$ is merely required to be disjoint from
$Q$, with no loss we may assume that the ends of
$\aaa_{\pm}$ on the circle $c = P \cap T$ are disjoint.  Say that the arcs
$\aaa_{\pm}$ {\em cross} if the ends of $\aaa_{+}$ and $\aaa_{-}$
alternate around $c$; that is, if the ends of
$\aaa_{+}$ lie on different arc components of $c - \aaa_{-}$.

\bigskip

{\bf Claim:} Some pair of arcs that satisfy the hypotheses for
$\aaa_{\pm}$ cross.

{\bf Proof of claim:} Assume, on the contrary, that no such pair of
arcs crosses.  Then among arcs of $Q \cap T_{\pm}$ satisfying the
conditions for $\aaa_{\pm}$ choose the pair whose ends are closest to
each other on the circle $c$.  The ends of $\aaa_{\pm}$ divide $c$
into four arcs, one of them, denoted $\bbb_{+}$, is bounded by the
ends of $\aaa_{+}$ and the other, denoted $\bbb_{-}$, by the ends of
$\aaa_{-}$.  Let $c_{\pm} = |Q \cap \bbb_{\pm}|$.

$T_{+} - \eta(\aaa_{+})$ is an annulus $A$; denote the boundary
component that contains $\bbb_{\pm}$ by $\bdd_{\pm} A$.  Then
$|\bdd_{-} A \cap Q| = c_{+}$ and $|\bdd_{-} A \cap Q| \geq c_{-}$. 
(The inequality reflects the fact that $Q$ may also intersect the two
intervals $c - \bbb_{\pm}$.)  No arc of $Q \cap A$ can have both ends
on $\bdd_{-} A$, else it would have been parallel to $\aaa_{+}$ in
$T_{+}$, and yet closer to $\aaa_{-}$.  We conclude that $c_{+} \geq
c_{-}$.  Arguing symmetrically on $T_{-} - \eta(\aaa_{-})$, we obtain
$c_{-} \geq c_{+} + 2$, the extra $2$ arising from the ends of
$\aaa_{+} \subset Q - c_{+}$.  The two inequalities conflict, a 
contradiction proving the claim. 

\bigskip

Following the claim, we assume that $\aaa_{\pm}$ cross.  Let $\rrr
\subset T$ be the circle obtained by banding the circle $c$ to itself
along the two arcs $\aaa_{\pm} \subset T_{\pm}$.  It is a single
circle because $\aaa_{\pm}$ cross.  Moreover, it's easy to see that
$\rrr$ is an essential circle in $T$ (there are essential curves in
$T$ on both sides of $\rrr$) and that $\rrr$ bounds disks both in $V$
and $W$.  So $\rrr$ is the intersection with $T$ of a reducing sphere
$R$.  Moreover, $R \cdot P = |\rrr \cap P| = 4$ and $R \cdot Q \leq
|\rrr \cap Q| = |c \cap Q| - 2 \leq P \cdot Q - 2$ since the ends of
$\aaa_{+}$ no longer count.
\end{proof}

\begin{prop} \label{prop:existsR} There is a reducing sphere $R$ so
that $P \cdot R = 4$ and $R \cdot Q < P \cdot Q$.
\end{prop}

\begin{proof} If there are two arcs of $(Q \cap T) - c$, one of slope
$0$ and one of slope $\infty$, one lying in $T_{+}$ and the other
lying in $T_{-}$, the result follows from Lemma \ref{lemma:cross}. 
Following Lemma \ref{lemma:outermost} we know that there arcs of slope
both $0$ and $\infty$.  Thus we are done unless both these arcs lie in
$T_{-}$ say, and each arc of $Q \cap T_{+}$ has finite, non-zero
slope.  Moreover, if all arcs of $Q \cap T_{+}$ have slope $1$ (or
slope $-1$) then a curve of slope $0$ in $T_{+}$ will be disjoint from
$Q \cap T_{+}$ and again we would be done by Lemma \ref{lemma:cross}. 
If $\sss, \ttt$ are slopes of arcs in $Q \cap T_{+}$, then, because
$|\Delta(\sss,\ttt)| \leq 1$, the inequality $0 < |\sss| < 1$ would
imply that $|\ttt| \leq 1$ and that $\sss$ and $\ttt$ have the same
sign.  Finally, a curve that has slope $\sss$, will have slope
$1/\sss$ if the roles of $V$ and $W$ are reversed.  Following these
considerations, we may as well restrict to the following case:

\begin{itemize}
    
    \item  Both slopes $0$ and $\infty$ arise among the arcs of $Q 
    \cap T_{-}$ and
    
    \item all arcs of $Q \cap T_{+}$ have slope $\sss$ with $0 < \sss
    \leq 1$ and not all have slope 1.
    
    \end{itemize}
    
Now consider a sphere $P^{+} \subset B_{+}$ that intersects the solid
torus $V \cap B_{+}$ in two meridian disks, and so intersects $W$ in
an annulus.  Again isotope the curve $Q \cap T$ so that it intersects
the two meridian circles $P^{+} \cap T$ minimally.  Any arc of $Q \cap
T_{+}$ must intersect $P^{+}$, else the arc would be of slope $0$.  In
particular, there is an essential non-separating disk $F \subset W$ so
that $\bdd F \subset T_{+}$ (i.  e. $\bdd F$ is a longitude of the
solid torus $V \cap B_{+}$) so that $F \cap P^{+}$ is a single
spanning arc of the annulus $P^{+} \cap W$ and so that the arc of
$\bdd F - P^{+}$ lying below $P^{+}$ (i.  e. in the pair of pants
component of $T_{+} - P^{+}$ adjacent to $c$) is disjoint
from $Q$.  See Figure \ref{fig:Pplus}

\bigskip

\begin{figure} [tbh]
\centering
\includegraphics[width=0.7\textwidth]{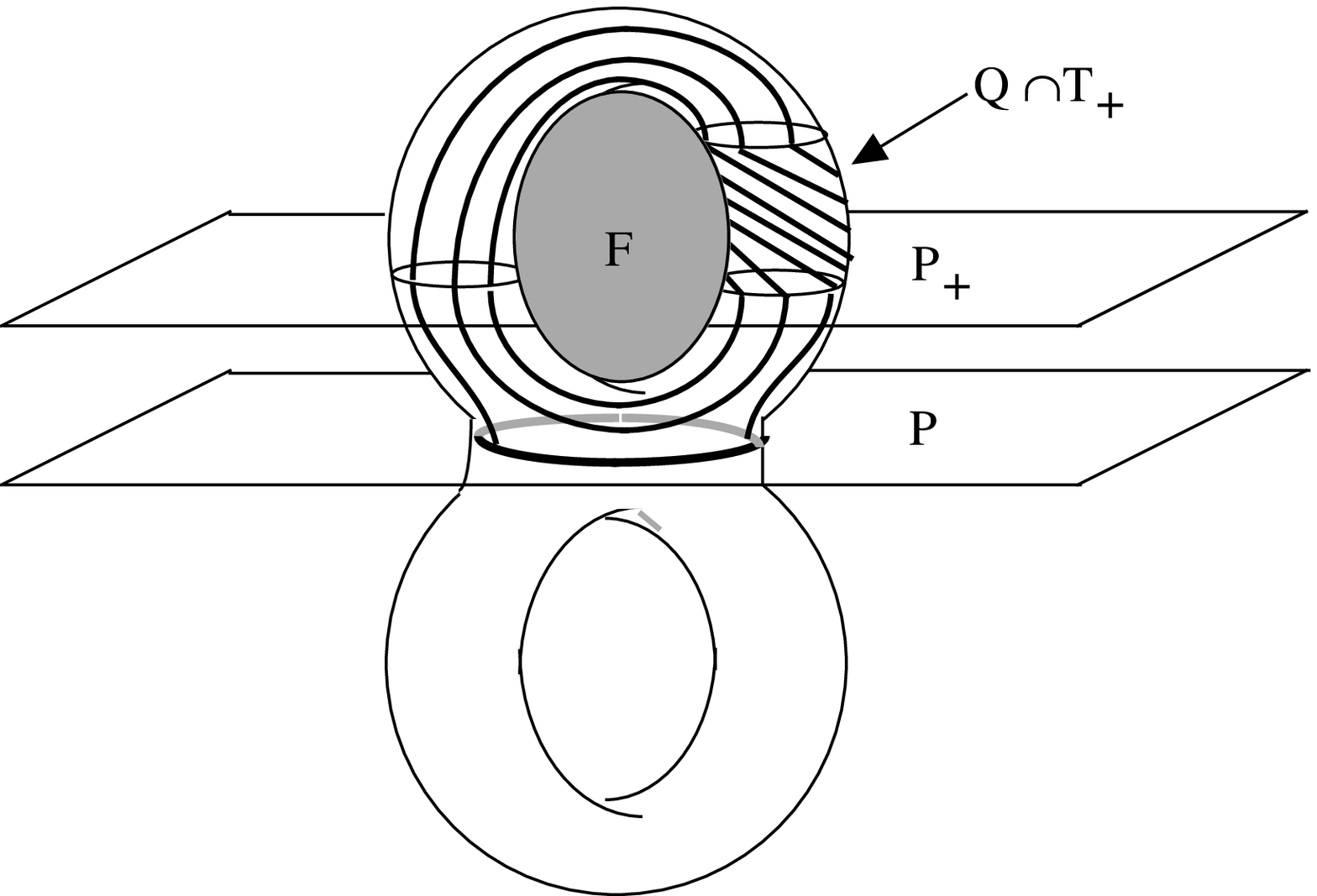}
\caption{} \label{fig:Pplus}
\end{figure}

We now examine outermost disks cut off from the disk $Q \cap W$ by the 
annulus $P^{+} \cap W$.  Let $E$ be any such disk.  Let $V^{\pm}$ be 
the closed components of $V - P^{+}$, with $V^{+}$ the $1$-handle  
lying above $P^{+}$ and $V_{-}$ the solid torus lying below $P^{+}$.  

{\bf Claim 1:} The outermost arc $\eee = \bdd E \cap P^{+}$ spans the
annulus $P^{+} \cap W$.

{\bf Proof of Claim 1:} This is obvious if $E$ lies above $P^{+}$,
since all arcs of $Q \cap T$ above $P^{+}$ span the $1$-handle
$V^{+}$.  If $E$ lies below $P^{+}$ the argument is a bit more subtle. 
Note that $V^{-}$ is a solid torus with two disks $d_{1}, d_{2}$ in
$\bdd V^{-}$ attached to $P^{+}$.  A simple counting argument (the
$d_{i}$ are parallel in $V^{+}$) shows that any arc of $Q \cap (\bdd
V^{-} - P^{+})$ that has both ends on the same disk $d_{i}$ is
essential in the torus $\bdd V^{-}$.  So an outermost disk $D \subset
V^{-}$ cut off from the disk $Q \cap V$ by the meridian disks $d_{i}$
must be a meridian disk of the solid torus $V^{-}$, with both ends on
$d_{1}$, say.  The same counting argument shows that some essential
arc in $Q \cap V^{-}$ must have both its ends on $d_{2}$ and so is a
meridinal arc for $V^{-}$ there as well.  If the ends of $\eee$ were
both on the same $d_{i}$, then $\eee$ would be a longitudinal arc
disjoint from the meridinal arc with ends at the other disk $d_{j}, j
\neq i$.  But a longtiudinal arc and a meridinal arc based at
different points must necessarily intersect.  Hence the ends of $\eee$
each lie on a different disk $d_{i}$, proving Claim 1.

\bigskip
    
{\bf Claim 2:} All the outermost disks cut off from $Q \cap W$ by
$P^{+}$ must lie on the same side of $P^{+}$.
    
    {\bf Proof of claim 2:} Suppose, on the contrary, that the
    outermost disks $E^{\pm}$ are cut off, with $E^{-}$ lying in the
    component of $S^{3} - P^{+}$ that lies
    below $P^{+}$ and $E_{+}$ lying in the component that lies above
    $P^{+}$.  Following Claim 1, both arcs $\eee^{\pm} = E^{\pm} cap 
    P^{+}$ span the annulus $P^{+} \cap W$.
    
Since the arc $E_{-} \cap T$ is disjoint from $\bdd F$ it follows from
a simple innermost disk, outermost arc argument, that all of $E_{-}$
can be made disjoint from $F$; in particular the spanning arcs
$\eee^{-}$ and $F \cap P^{+}$ are disjoint.  Since the spanning arc
$\eee^{+}$ is disjoint from the spanning arc $\eee^{-}$ which in turn
is disjoint from the spanning arc $F \cap P^{+}$, $\eee^{+}$ can be
isotoped off of $F \cap P^{+}$ without moving $\eee^{-}$.  (See Figure
\ref{fig:Eminus}.)  Then again an innermost disk, outermost arc
argument allows us to isotope all of $E_{-}$ off of $F$.  Now consider
any arc component $\ggg$ of $(Q \cap T_{+}) - P^{+}$.  If $\ggg$ lies
below $P^{+}$ then it is disjoint from $\bdd F$, by construction; if
$\ggg$ lies above $P^{+}$ then since it is disjoint from $E^{+}$, it
intersects $\bdd F$ at most once.  In particular, any arc of $Q \cap
T^{+}$ intersects a component of $P^{+} \cap T_{+}$ at least as often
as it intersects $\bdd F$, hence its slope has absolute value $\geq
1$.  This contradicts the second property itemized above, and so
proves claim 2).

\begin{figure} [tbh]
\centering
\includegraphics[width=0.7\textwidth]{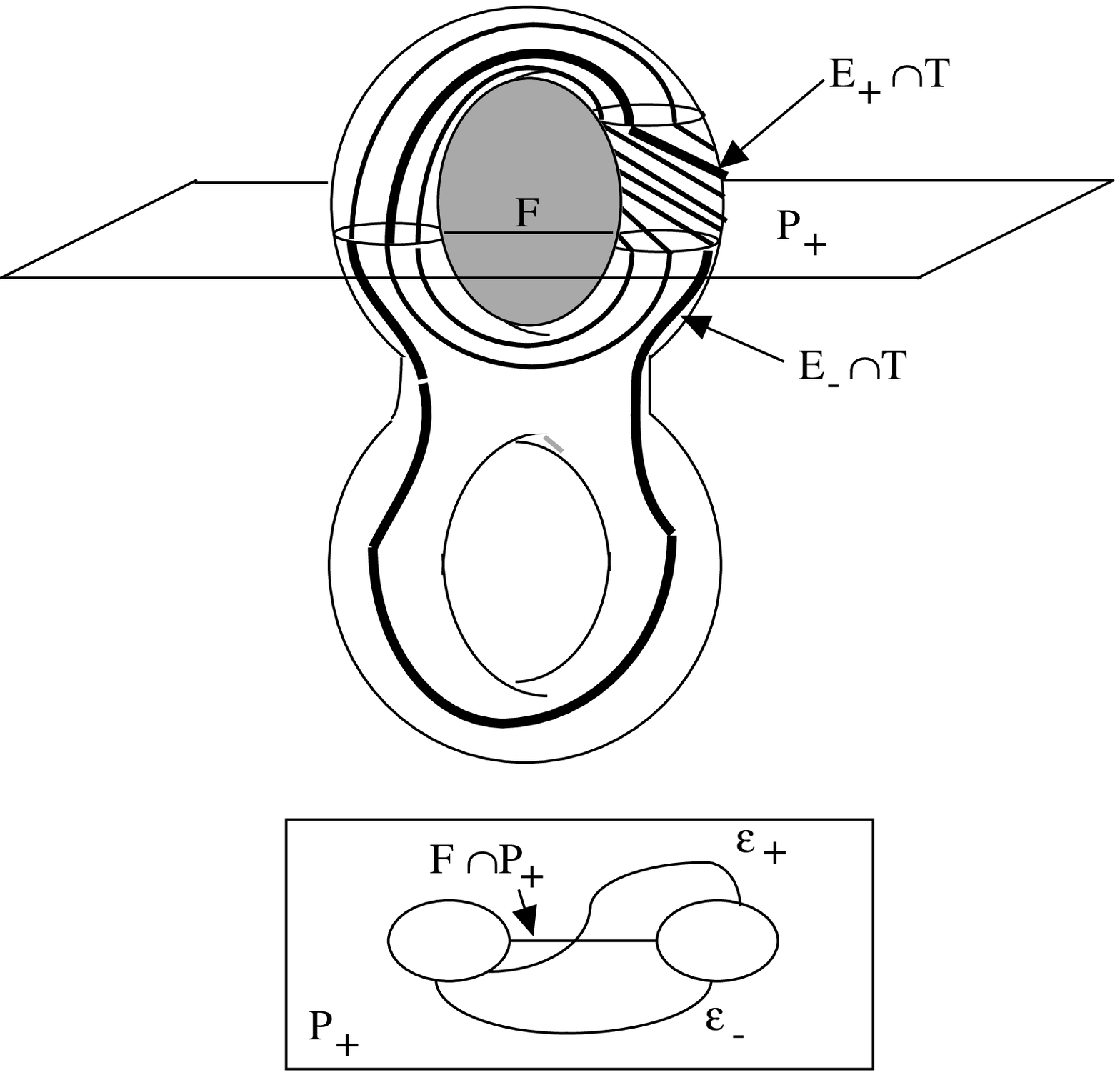}
\caption{} \label{fig:Eminus}
\end{figure}

\bigskip

{\bf Claim 3:} All the outermost disks cut off from $Q \cap W$ by
$P^{+}$ must lie above $P^{+}$.

{\bf Proof of claim 3:} Following claim 2) the alternative would be
that they all lie below (in $B_{-}$).  We show how this leads to a
contradiction.  Consider the disk $Q \cap W$ and how it is cut up by
the annulus $P^{+} \cap W$.  A standard innermost disk argument
ensures that all closed curves of intersection can be removed.  There
is at least one (disk) component $E_{0}$ of $(Q \cap W) - P^{+}$ that
is ``second outermost'', i.  e. it is adjacent to some $n \geq 2$
other components of $(Q \cap W) - P^{+}$, all but at most one of them
outermost.  See Figure \ref{fig:secondout}.  Since $E^{0}$ is adjacent to
an outermost component, all of which we are assuming lie below
$P^{+}$, $E_{0}$ must lie above $P^{+}$.  By Claim 1), all the
outermost arcs of intersection of $P^{+}$ with the disk $Q \cap W$
must span the annulus $W \cap P^{+}$, so it follows that each of the
$n$ arc components of $\bdd E_{0} \cap T$ spans the $1$-handle
$V^{+}$.  In particular, the union of the disk $E_{0}$ with the
punctured solid torus $P^{+} \cup V^{+}$ is the spine of a Lens space
in $S^{3}$, a contradiction proving Claim 3).

\begin{figure} [tbh]
\centering
\includegraphics[width=0.3\textwidth]{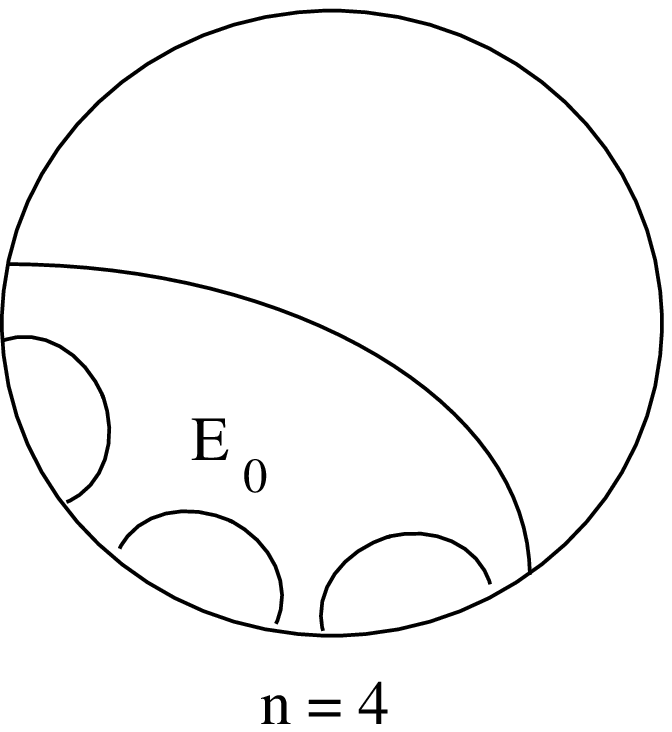}
\caption{} \label{fig:secondout}
\end{figure}

\begin{figure} [tbh]
\centering
\includegraphics[width=0.7\textwidth]{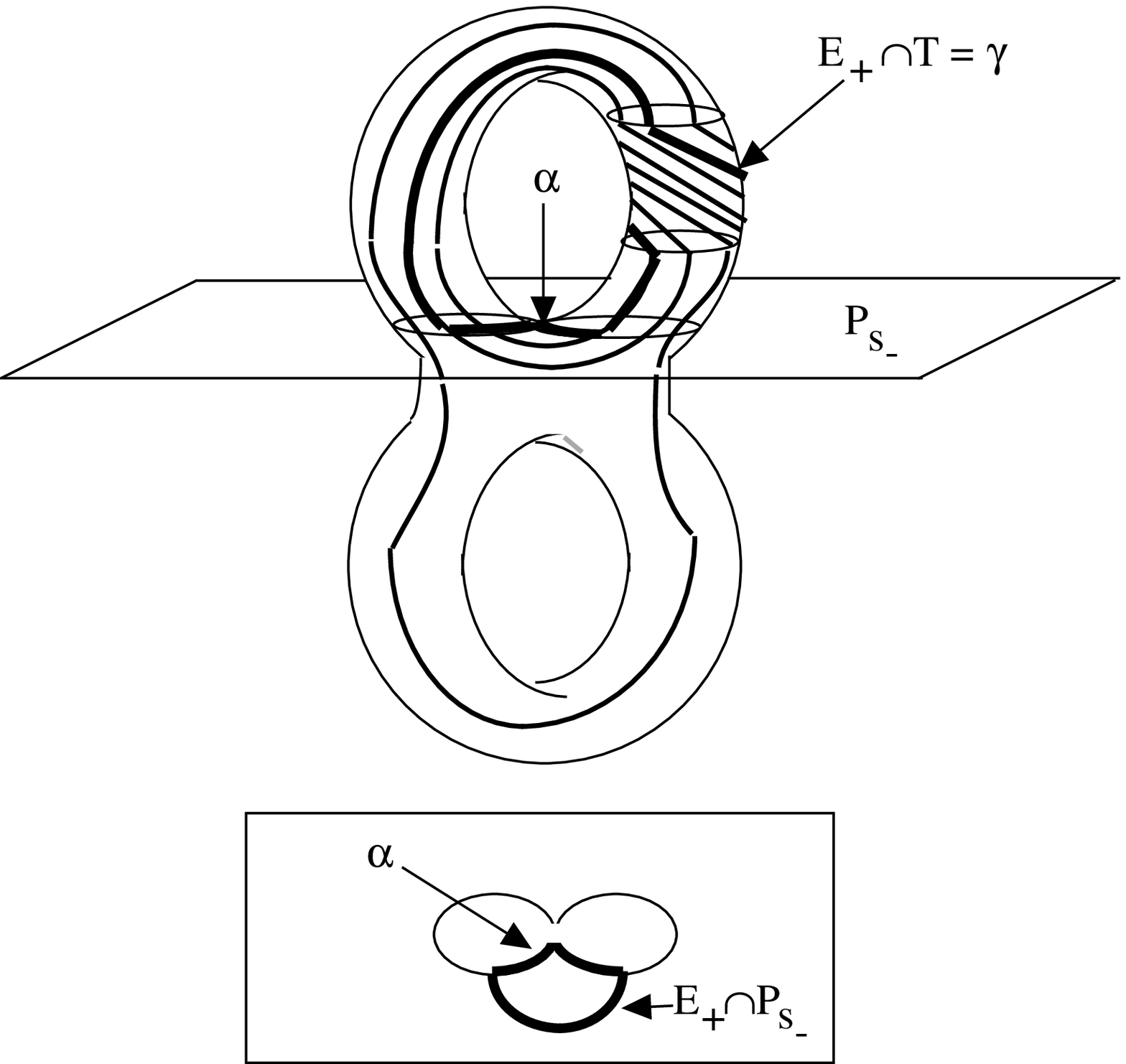}
\caption{} \label{fig:saddle}
\end{figure}

\bigskip

Following Claim 3), consider a sphere $P^{s}$ that passes through the
saddle point of $T_{+}$ that lies below $B^{+}$.  We can assume (see
Figure \ref{fig:saddle}) that $P^{s}$ intersects $Q$ transversally and
that every arc of $Q \cap T^{+}$ that lies above $P^{s}$ spans the
$1$-handle $V^{+}$.  According to claim 3) applied to a plane just
slightly higher than $P^{s}$, $P^{s}$ (and so also a plane $P^{s_{-}}$
lying just below $P^{s}$) cuts off a disk $E^{+}$ from $Q \cap W$ that
lies above the plane.  Let $\aaa \subset (P^{s_{-}} \cap T)$ be an arc
parallel in $P^{s_{-}} \cap W$ to the arc $E^{+} \cap P^{s_{-}}$, so
the union $\lll$ of $\aaa$ and the arc $\ggg = E^{+} \cap T$ is a
longitude lying above $P^{s_{-}}$ (indeed $\lll$ is a meridian of
$W$).  It's easy to isotope the ends of $\ggg$ closer together in
$\aaa$ until no arc of $(Q \cap T) - P^{s_{-}}$ lying above
$P^{s_{-}}$ has more than one end on $\lll$.  It then follows just as
in the proof of Claim 2) that any arc component of $Q \cap T_{+}$
intersects a meridian of $V^{+}$ at least as often as it intersects
$\lll$ and so has slope $\geq 1$, a contradiction that completes the
proof.  \end{proof}

\begin{cor} \label{cor:connected} The $2$-complex $\Ggg$ is connected.
    \end{cor}
    
     \begin{proof} Let $w$ be a fixed vertex of $\Ggg$, with
     associated reducing sphere $Q$.  Let $\Ggg_{0}$ be any component
     of $\Ggg$.  Choose a reducing sphere $P$ among those represented
     by vertices in $\Ggg_{0}$ so that $P \cdot Q$ is minimized. 
     Unless $P = Q$, Proposition \ref{prop:existsR} provides a
     reducing sphere $R$ which is represented by a vertex in
     $\Ggg_{0}$ (indeed one adjacent to the vertex representing $P$)
     but for which $R \cdot Q < P \cdot Q$.  From the contradiction we
     conclude then that indeed $P = Q$, so $w \in \Ggg_{0}$. 
     \end{proof}
     
Corollary \ref{cor:connected} is essentially \cite[Proposition
2.6]{ST}.  There we used Goeritz' theorem to prove the proposition;
here we have proven the proposition from first principles and now
observe that it proves Goeritz' theorem.

\section{A finite set of generators}

\begin{theorem} Suppose $\ddd \in \Hhh$ is any element with the
property that $P \cdot \ddd(P) = 4.$ Then the group $\Hhh$ is
generated by $\aaa, \bbb, \ggg, \ddd.$
\end{theorem}

\begin{proof} Choose any $h \in \Hhh$ and let $Q = h(P)$.  If $Q = P$
the result follows immediately from Lemma \ref{lemma:stabilize}. 
Otherwise, following Corollary \ref{cor:connected}, there is a
sequence of reducing spheres $P = P_{0}, P_{1}, \ldots, P_{n} = Q$ so
that $P_{i-1} \cdot P_{i} = 4, i = 1, \ldots, n$.  The proof will be
by induction on the length $n$ of this sequence -- the case $n = 1$
follows from Corollary \ref{cor:linktransitive}.  In particular, there
is a word $\ooo$ in the generators $\aaa, \bbb, \ggg, \ddd$ so that
$\ooo (P_{1}) = P$.  Apply $\ooo$ to every sphere in the shorter
sequence $P_{1}, \ldots, P_{n} = Q$ and obtain a sequence $P = \ooo
(P_{1}), \ooo (P_{2}), \ldots, \ooo(Q) = \ooo(h(P))$.  Then by
inductive hypothesis, $\ooo h$ is in the group generated by $\aaa,
\bbb, \ggg, \ddd$, hence so is $h$.
\end{proof}

There are several natural choices for $\ddd$.  For example, if we
think of $V$ as a ball with two $1$-handles attached, the two
$1$-handles separated by the reducing sphere $P$, then a slide of an
end of one of the $1$-handles over the other around a longitudinal
curve will suffice for $\ddd$.  This is the genus two version of
Powell's move $D_{\theta}$ (\cite[Figure 4]{Po}).  Another possibility
is to choose an order two element for $\ddd$, an element that is
conjugate in $\Hhh$ to $\ggg$: note from Figure \ref{fig:PandQandR}
that $Q \cdot \ggg(Q) = 4$.  

\begin{figure} [tbh]
\centering
\includegraphics[width=0.4\textwidth]{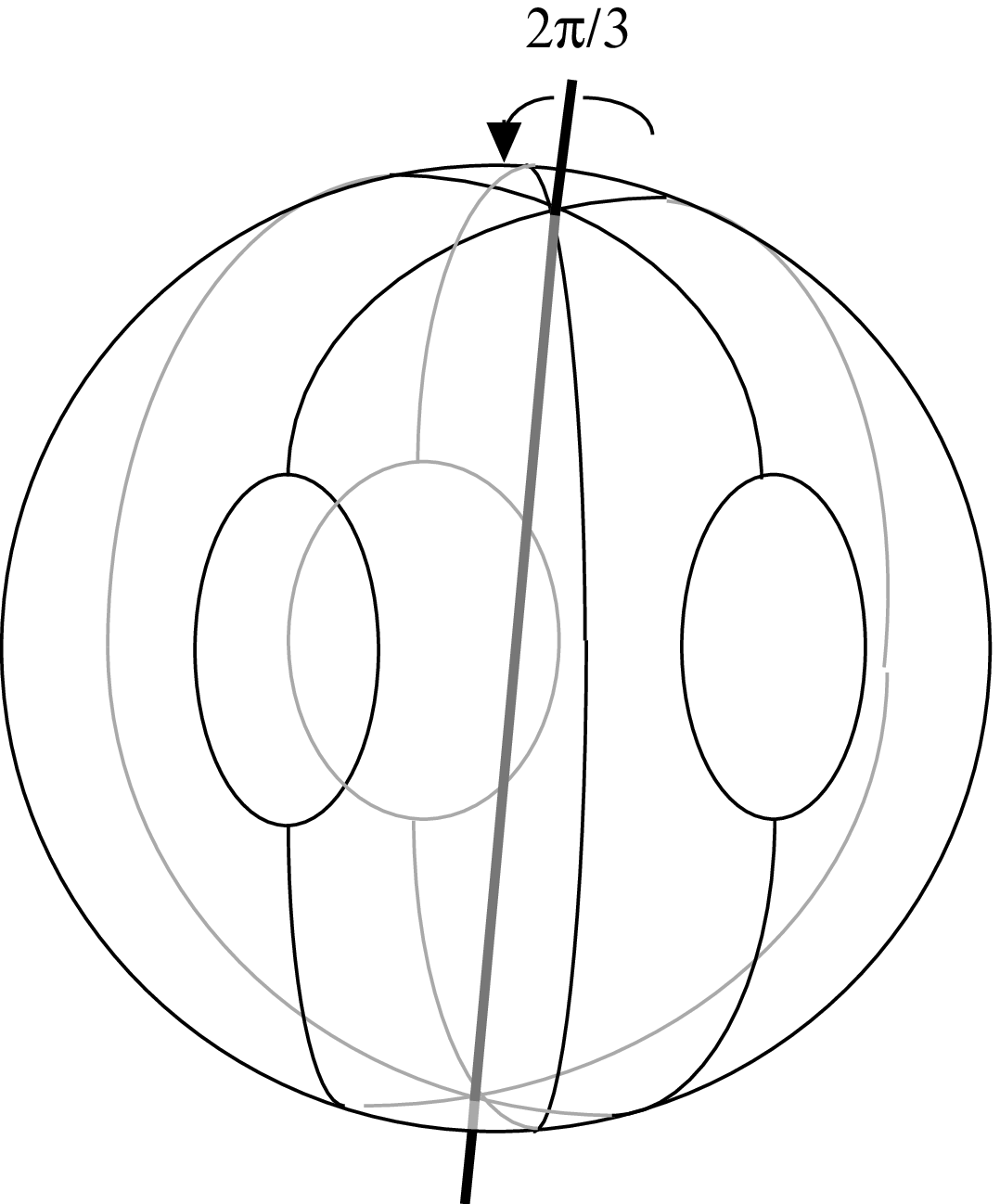}
\caption{} \label{fig:rotate}
\end{figure}

A bit more imaginative is the automorphism shown in Figure
\ref{fig:rotate} which is of order three and corresponds to rotating
one of the two-simplices of $\Ggg$ around its center.  The figure is
meant to evoke a more symmetric version of Figure \ref{fig:PandQandR}:
it depicts a thrice punctured sphere with three essential arcs, each
pair intersecting in two points.  Thicken the figure (i.e. cross with
an interval).  Then the thrice punctured sphere becomes a genus two
handlebody $V$ and each arc becomes a disk.  Each disk is the
intersection with $V$ of a reducing sphere, and the three reducing
spheres are represented by the corners of a single two-simplex $\sss$ in
$\Ggg$.  Rotation of the figure by $2\pi/3$ along the axis shown
cyclically permutes the three arcs, and so cyclically permutes the
three reducing spheres.  Hence it also rotates the corresponding 
$2$-simplex $\sss$ in $\Ggg$.  

\end{document}